\documentclass[12pt]{amsart}

\usepackage{amssymb, amsthm, hyperref}
\usepackage{enumerate, amsfonts, latexsym,epsfig, color}

\input xy
\xyoption{all}

\newtheorem {theorem}{Theorem} [section]
\newtheorem {lemma} [theorem] {Lemma}
\newtheorem {proposition} [theorem] {Proposition}
\newtheorem {example} [theorem] {Example}

\newtheorem {definition} [theorem] {Definition}

\newtheorem {remark} [theorem] {Remark}

\newtheorem {question}[theorem] {Question}

\newtheorem{thmA}{Theorem}

\def\co{\colon\thinspace}

\def\R {\mathbb R}
\def\N {\mathbb N}
\def\ssm {\smallsetminus}
\def\mc {\mathcal}

\def\Z {\mathbb{Z}}
\def\curlyR {\mathcal{R}}

\def\Hom {\mathrm{Hom}}
\def\Aut {\mathrm{Aut}}

\def\IHom {\text{IHom}_\Gamma(H_\Sigma, \Gamma ; \Lambda)}
\def\SK {\underrightarrow{\text{Ker}}\, }
\def\O {\mc{O}(H_\Sigma, \Gamma; \Lambda)}

\def\IHomL {\text{IHom}_\Gamma(\hat{L}, \Gamma ; \Lambda)}
\def\IHomG {\text{IHom}_\Gamma(G, \Gamma ; \Lambda)}
\def\OL {\mc{O}(\hat{L}, \Gamma; \Lambda)}
\def\OG {\mc{O}(G, \Gamma; \Lambda)}

\begin{document}

\title[Conjugacy classes of solutions]
{Conjugacy classes of solutions to equations and inequations over hyperbolic groups}

\author[Daniel Groves]{Daniel Groves}
\address{Daniel Groves\\
MSCS UIC 322 SEO, \textsc{M/C} 249\\
851 S. Morgan St.\\
Chicago, IL 60607-7045, USA}
\email{groves@math.uic.edu}

\author[Henry Wilton]{Henry Wilton}
\address{Henry Wilton\\
Dept of Mathematics\\
1 University Station C1200\\
Austin, TX 78712, USA}
\email{henry.wilton@math.utexas.edu}

\date{\today}

\thanks{The first author's work was supported in part by NSF Grant DMS-0504251.  Much of the work was undertaken while the first author was a Taussky--Todd Instructor at Caltech.}

\begin{abstract}
We study conjugacy classes of solutions to systems of equations and inequations over torsion-free hyperbolic groups, and describe an algorithm to recognize whether or not there are finitely many conjugacy classes of solutions to such a system.  The class of immutable subgroups of hyperbolic groups is introduced, which is fundamental to the study of equations in this context. We apply our results to enumerate the immutable subgroups of a torsion-free hyperbolic group.
\end{abstract}

\maketitle

\section{Introduction}

This paper is concerned with decision problems relating to systems of equations and inequations over torsion-free hyperbolic groups, a subject with a long history.  Makanin \cite{Makanin} proved that it is possible to decide whether a system of equations over a free group has a solution, and in \cite{Makanin2} he extended this result to systems of equations and inequations.  Razborov \cite{Razborov, Razborov:thesis} described an algorithm whose output gives a parametrization of the entire set of solutions to a system of equations over a free group. His description was further refined by Kharlampovich and Miasnikov \cite{KM}.

Over torsion-free hyperbolic groups, Rips and Sela \cite{RS} proved that it is possible to decide if a system of equations has a solution.  The corresponding result for systems of equations and inequations was proved independently by Sela \cite{Sela:hyp} and Dahmani \cite{Dah_eq}.

The study of equations over groups has recently received a lot of attention due to the solutions to Tarski's problem by Sela \cite{Sela:last} and also by Kharlampovich and Miasnikov \cite{KM:Final}.  Our approach is inspired by Sela's results, and we use some of his structure theory.

\begin{thmA} \label{t:FinSol}
There is an algorithm that takes as input a finite presentation for a torsion-free hyperbolic group $\Gamma$ and a finite system of equations and inequations with coefficients in $\Gamma$ and determines whether or not the system has only finitely many solutions in $\Gamma$.

If there are  finitely many solutions, the algorithm provides a list of all the solutions.
\end{thmA}

In the case where $\Gamma$ is a free group and the system has only equations, this theorem is immediate from Razborov's description of the set of solutions.  Also, if $\Gamma$ is free and there are inequations, it is not difficult to deduce Theorem \ref{t:FinSol} from Razborov's results.  However, for general torsion-free hyperbolic $\Gamma$, more work is required.

When the system of equations and inequations has no coefficients, $\Gamma$ naturally acts by conjugation on the set of solutions.  Therefore, if $\Gamma$ is non-elementary and there is at least one solution then there are infinitely many.  Thus, in the case of coefficient-free equations and inequations, Theorem \ref{t:FinSol} follows immediately from the above-mentioned result of Sela and Dahmani.

It is therefore natural to consider the orbits of solutions under the action by conjugation.  For example, in Sela's definition of limit groups, he takes a minimal solution in a conjugacy class in order to ensure that the limiting action on an $\R$-tree has no global fixed point.  Conjugacy classes of solutions are thus of fundamental concern for Sela's approach to the Tarski problem (although this is only the very beginning of the story).

Our next result concerns conjugacy classes of solutions.

\begin{thmA} \label{t:FinConjSol}
There is an algorithm that takes as input a finite presentation for a torsion-free hyperbolic group $\Gamma$ and a finite (coefficient-free) system of equations and inequations and determines whether or not the system has only finitely many conjugacy classes of solutions in $\Gamma$.

If there are finitely many conjugacy classes of solutions, the algorithm provides a list containing exactly one solution in each conjugacy class.
\end{thmA}

\begin{remark}
If $\Gamma$ is free, when there is one solution (and no coefficients) there are infinitely many conjugacy classes of solutions.  However, this may not be the case for arbitrary hyperbolic $\Gamma$.  See Example \ref{e:Finitely many conjugacy classes of solutions}.
\end{remark}

Although we believe it to be of independent interest, one of our motivations for proving Theorem \ref{t:FinConjSol} is its use in the study of $\Gamma$-limit groups (see Section \ref{s:Limit groups}), where $\Gamma$ remains a torsion-free hyperbolic group.  Unlike the case of limit groups (when $\Gamma$ is free), $\Gamma$-limit groups may not be finitely presentable.  In fact, every $\Gamma$-limit group is finitely presentable if and only if $\Gamma$ is coherent.

The basic building blocks of $\Gamma$-limit groups are free groups and finitely generated subgroups of $\Gamma$ that admit only finitely many conjugacy classes of embeddings into $\Gamma$.  We call such subgroups {\em immutable}.  Any study of $\Gamma$-limit groups, or even the subgroups of $\Gamma$, must take into account immutable subgroups.  They are one of the key technical differences between ordinary limit groups over free groups and $\Gamma$-limit groups (see Remark \ref{r:Immutable vs Strict}).

The expert reader may be interested to learn that, by the results of \cite{Sela:hyp}, there is an analogue over torsion-free hyperbolic groups of the analysis lattice that appears in \cite{Sela:IHES}.  Immutable subgroups appear at the lowest level of the analysis lattice.  This analogy will be made explicit in forthcoming work \cite{GW3}.

We will see in Section \ref{s:Immutable} that there is a natural relationship between immutable subgroups and systems of equations and inequations with finitely many conjugacy classes of solutions (see Example \ref{e:Finitely many conjugacy classes of solutions}).   We exploit Theorem \ref{t:FinConjSol} to enumerate the immutable subgroups of $\Gamma$.

\begin{thmA} \label{t:ListImm}
Let $\Gamma$ be a torsion-free hyperbolic group.  The set of immutable subgroups of $\Gamma$ is recursively enumerable, in the following sense: there exists a Turing machine that lists all finite subsets of $\Gamma$ that generate immutable subgroups.
\end{thmA}

This result is perhaps surprising because the finitely generated subgroups of (torsion-free) hyperbolic groups are notoriously badly behaved when it comes to decision problems.  However, one of the themes of our work is that when it comes to equations and inequations, most of these difficulties can be surmounted.  In \cite{GW3}, we will combine Theorem \ref{t:ListImm} with the techniques of \cite{GW1} to recursively enumerate all $\Gamma$-limit groups.

\begin{remark}
All of the algorithms we provide are uniform, in the sense that they can take as input a finite presentation for $\Gamma$ as well as the equations and inequations.  Thus, they work for all torsion-free hyperbolic groups in a uniform manner.
\end{remark}

An outline of this paper is as follows.  In Section \ref{s:EqHom} we record the obvious rephrasing of solutions to equations and inequations over $\Gamma$ in terms of homomorphisms to $\Gamma$.  Section \ref{s:Limit groups} contains some basic definitions and facts about $\Gamma$-limit groups.  In Section \ref{s:Pseudo} we introduce the notion of {\em $\Gamma$-approximations}, which are a key tool to circumvent the fact that $\Gamma$-limit groups may not be finitely presentable.  In Section \ref{s:Infinite} we present the halves of the algorithms required for Theorems \ref{t:FinSol} and \ref{t:FinConjSol} which recognize when the appropriate sets are infinite, whereas in Section \ref{s:Finite} we provide the algorithms which recognize when they are finite.  Finally, in Section \ref{s:Immutable} we discuss immutable subgroups of $\Gamma$ and prove Theorem \ref{t:ListImm}.  We end with some questions about immutable subgroups.

\section{Equations and homomorphisms} \label{s:EqHom}

The class of (word-)hyperbolic groups is well known---for an introduction, see for instance \cite[Part III]{Bridson-Haefliger}. Throughout this paper $\Gamma$ will be a fixed non-elementary torsion-free hyperbolic group.  (The case in which $\Gamma$ is elementary is interesting, but classical.)  Fix a finite generating set $\mc{A} = \{ a_1 , \ldots , a_k \}$ for $\Gamma$.  Let $F(\mc{A})$ be the free group on $\mc{A}$, and $\theta \co F(\mc{A}) \to \Gamma$ the canonical epimorphism.  We will blur the distinction between the $a_i$ as elements of $F(\mc{A})$ and as elements of $\Gamma$.

Let $X = \{ x_1, \ldots , x_n \}$ be a finite set of variables, and $F(X)$ the free group on $X$.

A {\em system of equations} is a collection $\Sigma \subseteq F(X) \ast F(\mc{A})$.  Write $\Sigma = \{ \sigma_i(\underline{x}, \underline{a}) \}_{i \in I}$.  A {\em solution} to $\Sigma$ is a collection $\{ \gamma_k \}_{k = 1}^n \subseteq \Gamma$ so that for each $i \in I$ we have $\sigma_i(\underline{\gamma}, \underline{a}) =_\Gamma 1$.  In other words, the result of substituting $\gamma_j$ for $x_j$ in each $\sigma_i$  (and interpreting the $a_k$ as the fixed generators for $\Gamma$) yields the trivial   element of $\Gamma$.

A {\em system of equations and inequations} is a pair of subsets $\Sigma, \Lambda \subseteq F(X) \ast F(\mc{A})$.  Suppose that $\Lambda = \{ \lambda_j(\underline{x},\underline{a}) \}_{j \in J}$. A {\em solution} to $\Sigma = 1$, $\Lambda \ne 1$ is a subset $\underline{\gamma} \subseteq \Gamma$ so that for each $i \in I$ we have $\sigma_i(\underline{\gamma},\underline{a}) = 1$ and for each $j \in J$ we have $\lambda_j(\underline{\gamma},\underline{a}) \ne 1$.

The letters $a_i$ (and $a_i^{-1}$) appearing in elements of $\Sigma$ and $\Lambda$ are called {\em coefficients}.  We say that $\Sigma$ and $\Lambda$ are {\em missing coefficients} if they are subsets of the canonical copy of $F(X)$ in $F(X) \ast F(\mc{A})$.  In case $\Sigma$ and $\Lambda$ are missing coefficients there are analogous subsets (still denoted $\Sigma$ and $\Lambda$) of $F(X)$, which we consider to be a {\em coefficient-free system of equations and inequations} (also called a system {\em without coefficients}).  Therefore, we make a distinction between systems of equations and inequations which happen not to contain any coefficients, and those which are coefficient-free.  This distinction is made because in the coefficient-free setting $\Gamma$ acts by conjugation on the set of solutions, whereas it does not when there are coefficients.

Following Razborov \cite{Razborov}, it is natural to interpret solutions to systems of equations and inequations as homomorphisms to $\Gamma$.  
If $\Sigma$ is coefficient-free, define the group $H_\Sigma$ by the presentation 
\begin{equation} \label{eq:Coefficientless Group}
H_\Sigma = \langle X \mid \sigma_i(\underline{x}) = 1, i \in I \rangle.
\end{equation}

If $\Sigma$ has coefficients, define the group $H_\Sigma$ by the presentation
\begin{equation} \label{eq:Coefficient Group} 
\langle X, \mc{A} \mid  \sigma_i(\underline{x},\underline{a}) = 1,{i \in I} \rangle.    
\end{equation}
In this case, the {\em coefficient subgroup} of $H_\Sigma$ is the subgroup generated by $\mc{A} \subseteq H_\Sigma$.

When $\Sigma$ has coefficients, the set of solutions to $\Sigma$ is naturally in one-to-one correspondence with the set of homomorphisms:
\[  \Hom_{\Gamma}(H_\Sigma,\Gamma) = \{ h \co H_\Sigma \to \Gamma \mid h(a_i) = a_i \} .   \]
If $\Sigma$ is coefficient-free the set of solutions is parameterized by the full set of homomorphisms from $H_\Sigma$ to $\Gamma$, denoted $\Hom(H_\Sigma,\Gamma)$.

The set of solutions to $\Sigma = 1$, $\Lambda \ne 1$ is naturally in one-to-one correspondence with a smaller set of homomorphisms:
\[  \IHom = \{ h \co H_\Sigma \to \Gamma \mid h(a_i) = a_i, h(\lambda_j) \ne 1 \}   .   \] (In the coefficient-free case, there is an obvious analogous definition for $\mathrm{IHom}(H_\Sigma,\Gamma;\Lambda)$.)

Under this simple reformulation, if $\Sigma$ and $\Lambda$ are finite, Theorem \ref{t:FinSol} asserts that there is an algorithm that decides if $\IHom$ is finite.

If $\Sigma$ and $\Lambda$ are coefficient-free, we can post-compose an element of $\mathrm{IHom}(H_\Sigma,\Gamma;\Lambda)$ by an element of $\text{Aut}(\Gamma)$, and we will get another element of $\mathrm{IHom}(H_\Sigma,\Gamma;\Lambda)$.  In particular, $\Gamma$ acts on $\mathrm{IHom}(H_\Sigma,\Gamma;\Lambda)$ by conjugation.

We denote the set $\mathrm{IHom}(H_\Sigma,\Gamma;\Lambda) / \Gamma$ by $\O$.  Theorem \ref{t:FinConjSol} asserts that there is an algorithm that recognizes whether or not $\O$ is finite.

\section{$\Gamma$-limit groups}\label{s:Limit groups}

One of our key tools is the theory of {\em $\Gamma$-limit groups}, as introduced and studied when $\Gamma$ is torsion-free hyperbolic by Sela in \cite{Sela:hyp}.  One of the main difficulties in studying $\Gamma$-limit groups as opposed to ordinary limit groups (when $\Gamma$ is free) is that if $\Gamma$ is not coherent then $\Gamma$-limit groups will not all be finitely presented.

\begin{remark}
The Rips construction \cite{Rips} can be used to construct many torsion-free hyperbolic groups which have finitely generated subgroups that are not finitely presented.
\end{remark}

Recall the following two definitions, due to Bestvina and Feighn in the case of free groups.

\begin{definition} \cite[Definition 1.5]{BFSela}
Let $G$ be a finitely generated group, and let $\{ h_n : G \to \Gamma \}$ be a sequence of homomorphisms.  The {\em stable kernel} of $\{ h_n \}$, denoted $\SK(h_n)$, is the set of all elements $g \in G$ so that $g \in \text{ker}(h_n)$ for all but finitely many $n$.

The sequence $\{ h_n \}$ is {\em stable} if for all $g \in G$, either $g \in \SK(h_n)$ or, for all but finitely many $n$, $g \not\in \text{ker}(h_n)$.
\end{definition}

\begin{definition} \cite[Definition 1.5]{BFSela} \label{DefAlgLimit}
A {\em $\Gamma$-limit group} is a group of the form $L = H_\Sigma/\SK(h_n)$ where $\Sigma$ is a coefficient-free system of equations, $H_\Sigma$ is as in
the previous section, and $\{ h_n \} \subseteq \Hom(H_\Sigma,\Gamma)$ is a stable sequence.

A {\em restricted $\Gamma$-limit group} is a group of the form $L = H_\Sigma / \SK(h_n)$ where $\Sigma$ is a system of equations {\em with} coefficients, and
$\{ h_n \} \subseteq \Hom_{\Gamma}(H_\Sigma,\Gamma)$ is a stable sequence.  The image of the coefficient subgroup of $H_\Sigma$ in $L$
is called the {\em coefficient subgroup} of $L$.
\end{definition}

It should be noted that restricted $\Gamma$-limit groups are, in particular, $\Gamma$-limit groups.  Almost all of the theory of restricted $\Gamma$-limit
groups is identical to that of $\Gamma$-limit groups.  Often in this paper we will omit the adjective `restricted'; this should not cause
the reader any confusion.

\begin{remark}
By taking a constant sequence of homomorphisms, we see that all finitely generated subgroups of $\Gamma$ are $\Gamma$-limit groups.
\end{remark}

\begin{definition}
A $\Gamma$-limit group is called {\em strict} if it can be represented as $H_\Sigma/\SK(h_n)$ where $\{ h_n \co H_\Sigma \to \Gamma \}$ are pairwise non-conjugate.

A restricted $\Gamma$-limit group is {\em strict} if it can be represented as $H_\Sigma/\SK(h_n)$ where $\{ h_n \co H_\Sigma \to \Gamma \}$ are pairwise distinct.
\end{definition}

We refer the reader to \cite{Sela:hyp} for the basic theory of $\Gamma$-limit groups where $\Gamma$ is a torsion-free hyperbolic group.  However, we do emphasize two features.  Recall that if $H$ and $G$ are groups then
we say that $G$ is {\em fully residually $H$} if for every finite subset $A \subseteq G$ there is a homomorphism
$h_A : G \to H$ which is injective on $A$.

\begin{proposition} \cite[Proposition 1.18]{Sela:hyp}
Let $\Gamma$ be a non-elementary torsion-free hyperbolic group.  A finitely generated group $G$ is fully
residually $\Gamma$ if and only if it is a $\Gamma$-limit group.
\end{proposition}
The above property will be used several times in this paper without comment.

  A strict $\Gamma$-limit group $L$ is naturally equipped with a faithful action on an $\R$-tree.  If $L$ is a strict restricted $\Gamma$-limit group, the coefficient subgroup of $L$ acts elliptically on this tree.  In order to apply the Rips Machine \cite{Sela:Inventiones,Guirardel} and deduce that $L$ splits (relative to the coefficient subgroup), we need
to know that the action has no global fixed point.  This is the content of the following standard proposition.  In the coefficient-free case, this is proved in
\cite{Paulin}, \cite{Bestvina} or \cite{BS}, for example.  The case with coefficients is well-known to the experts, but as far as we are aware there is no proof in the
literature.  We take this opportunity to record a proof, which was communicated to the first author by Sela.

\begin{proposition}\label{p:Action on a tree}
Any strict $\Gamma$-limit group $L$ acts without a global fixed point on an $\R$-tree. If $L$ is restricted the action may be chosen so that the coefficient subgroup fixes a point.
\end{proposition}
\begin{proof} [Proof (Case with coefficients)]
We start by outlining the standard construction of the $\R$-tree as a pointed Gromov-Hausdorff limit;  see
\cite{BS} or \cite{Bestvinahandbook}, for instance, for details.
Fix a word metric $d$ on $\Gamma$.  Let $\{ h_n \} \subseteq \Hom_\Gamma (H_\Sigma,\Gamma)$ be a stable sequence of pairwise distinct homomorphisms.  For each $n$ define the scaling factors $\delta_n$ and the basepoint $*_n=1\in \Gamma$.  For each $n$ we have an action of $H_\Sigma$ on $\Gamma$, equipped with the rescaled word metric $d/\delta_n$.  A subsequence of this sequence of actions converges (in the pointed Gromov--Hausdorff topology) to an action of $H_\Sigma$ on a Gromov-hyperbolic space $T$ equipped with a basepoint $*$.  Because the sequence $\{h_n\}$ is pairwise distinct, the scaling factors $\delta_n$ tend to infinity, and it follows that $T$ is an $\R$-tree.   The kernel of the action of $H_\Sigma$ on $T$ clearly contains $\SK(h_n)$, so the action of $H_\Sigma$ descends to an action of $L$ on $T$.

For each element $g$ of the coefficient subgroup of $H_\Sigma$, $h_n(g)=g$ for each $n$ and so $d(1,h_n(g))/\delta_n$  tends to $0$.  It follows that $g*=*$, so the coefficent subgroups of $H_\Sigma$ and hence $L$  act elliptically, as required. 

It remains to prove that the action of $L$ on $T$ has no global fixed point.  This is the part of the proof that differs significantly from the coefficient-free case.  The problem is that, because of the coefficient subgroup, we could not conjugate the homomorphisms $h_n$ to make $*_n$ centrally located.  Instead, we argue as follows.

The proof of \cite[Lemma 1.3(vi)]{Sela:hyp} applies to show that segment stabilizers of the action of $L$ on $T$ are abelian.  The coefficient subgroup of $L$ is clearly isomorphic to $\Gamma$ so is non-abelian.  Therefore, the basepoint $*$ is the unique point fixed by the coefficient subgroup.  On the other hand, by construction some generator $x_i$ moves $*$ distance $1$, so there is no global fixed point.  
\end{proof}

We are particularly interested in cyclic splittings that give rise to infinitely many outer automorphisms.

\begin{definition}
A cyclic splitting is \emph{essential} if the edge group is of infinite index in any vertex group.
\end{definition}

The point of Theorem \ref{t:CyclicSplit} is that one can choose the splitting assured by the Rips Machine to be essential.

\begin{theorem}[Sela, \cite{Sela:hyp}] \label{t:CyclicSplit}
Let $\Gamma$ be a non-elementary torsion-free hyperbolic group, and let $L$ be a non-abelian, freely indecomposable, strict $\Gamma$-limit group.  Then $L$ admits an essential cyclic splitting.

If, furthermore, $L$ is a restricted $\Gamma$-limit group then the splitting may be chosen so that the coefficient subgroup is elliptic.
\end{theorem}

Sela needs Theorem \ref{t:CyclicSplit} in \cite[Section 2]{Sela:hyp} (and claims that it is true on page 8).  However, he does not prove it in \cite{Sela:hyp}.  In \cite[Theorem 3.2]{Sela:IHES} it is proved when $\Gamma$ is free.  There are only a couple of extra subtleties when $\Gamma$ is not free, and essentially the same proof works (in particular, we know that Sela has a proof of Theorem \ref{t:CyclicSplit}).  However, we feel uncomfortable simply asserting that this strategy works so we include our own proof here.  The proof we present is inspired by a remark towards the beginning of the proof of \cite[Theorem 3.2]{Sela:hyp}, that Theorem \ref{t:CyclicSplit} is straightforward once one knows that all abelian subgroups of $L$ are finitely generated.

\begin{proof}[Proof of Theorem \ref{t:CyclicSplit}]
Note that $\Gamma$ is a toral relatively hyperbolic group.  Therefore we may apply the results of \cite{MR-RH} to $L$.  The advantage of this approach is that Theorem \ref{t:CyclicSplit} is false in this context (see \cite[Remark 5.1]{MR-RH}).  Therefore, the proof of \cite[Corollary 5.12]{MR-RH} does not rely on the result we are currently trying to prove, so there is no risk of circularity.  The content of \cite[Corollary 5.12]{MR-RH} is that abelian subgroups of $L$ are finitely generated.

Consider the abelian JSJ decomposition $\Delta$ of $L$ and collapse all edges that give rise to non-essential cyclic splittings of $L$, yielding a graph of groups $\Delta'$ for $L$.  As pointed out in the proof of \cite[Theorem 3.2]{Sela:IHES}, the action on an $\R$-tree that arises from Proposition \ref{p:Action on a tree} leads (via \cite[Theorem 1.5]{Sela:hyp}) to either an essential cyclic splitting or a non-trivial abelian splitting of $L$.  Therefore, we may assume that $\Delta'$ is non-trivial.

Suppose that $\Delta'$ contains only non-cyclic edge groups, and let $V$ be a non-abelian vertex group of this decomposition.  Then $V$ contains a non-cyclic abelian subgroup, so cannot be a subgroup of $\Gamma$.  Hence $V$ is a strict $\Gamma$-limit group, and so admits a non-trivial abelian splitting $\Theta$ in which all non-cyclic abelian subgroups are elliptic, by \cite[Lemma 1.7]{Sela:hyp}.  As noted above, $\Theta$ cannot be an inessential cyclic splitting.  The decomposition $\Theta$ of $V$ can be used to refine $\Delta'$.  If $\Theta$ has a trivial edge group, then we obtain a contradiction to the fact that $L$ is freely indecomposable.  Otherwise, $V$  is obtained from a 
rigid vertex group $W$ of $\Delta$ by adjoining finitely many roots and each such root is elliptic in every abelian splitting of $L$.  Therefore, $\Theta$ can be used to refine
$\Delta$, and we obtain an abelian splitting of $L$ in which $W$ is not elliptic, contradicting the canonical properties of $\Delta$.

In the case with coefficients, the proof is almost identical.  The only modification is that one works with the {\em restricted} abelian JSJ decomposition, rather than the full abelian JSJ decomposition.  The restricted decomposition is obtained by only considering those splittings in which the coefficient subgroup is elliptic.
\end{proof}

\section{$\Gamma$-approximations and splittings} \label{s:Pseudo}

We have already noted that $\Gamma$-limit groups need not be finitely presented.  Implementing algorithms involving groups that are not finitely presented is obviously more difficult than when the groups are finitely presented.  We will circumvent this difficulty by working with finitely presented approximations to our $\Gamma$-limit groups.

\begin{remark}
Any homomorphism $\phi \co G \to H$ induces a map \[\phi^* \co \Hom(H,\Gamma) \to \Hom(G,\Gamma)\] via $\phi^*(h) = h \circ \phi$.  If $\phi$ is an epimorphism then $\phi^*$ is an injection.
\end{remark}

\begin{definition}
Suppose that $G$ is a finitely generated group.    A {\em $\Gamma$-approximation to $G$} is a finitely presented
group $\hat{G}$ together with an epimorphism $\eta \co \hat{G} \to G$ so that $\eta^*$ is a bijection.
\end{definition}
Often the map $\eta$ will be implicit, and we will simply refer to $\hat{G}$ as a $\Gamma$-approximation to $G$.

We also need a version of $\Gamma$-approximations for equations with coefficients.

\begin{definition}
Let $\Sigma, \Sigma_0$ be systems of equations with coefficients, and $H_\Sigma, H_{\Sigma_0}$ the associated groups.  If $\Sigma_0 \subseteq \Sigma$ then there is a natural epimorphism $\eta \co H_{\Sigma_0} \to H_\Sigma$
which in turn induces an injection
\[	\eta^* \co \Hom_\Gamma(H_\Sigma,\Gamma) \hookrightarrow \Hom_{\Gamma}(H_{\Sigma_0},\Gamma)	.	\]
A {\em $\Gamma$-approximation to $H_\Sigma$} is a pair $(H_{\Sigma_0},\eta)$ where $\Sigma_0$ is a finite
system of equations and $\eta^*$ is a bijection.
\end{definition}

Saying that $\eta^*$ is a bijection is equivalent to saying that every homomorphism from $\hat{G}$ to $\Gamma$ factors through $\eta$.  To prove the existence of approximations, we appeal to the following result from \cite{Sela:hyp}, which states that torsion-free hyperbolic groups are {\em equationally Noetherian}.

\begin{theorem} \cite[Theorem 1.22]{Sela:hyp} \label{t:Noether}
Let $\Gamma$ be a torsion-free hyperbolic group.  Then every infinite system of equations in finitely many variables over $\Gamma$ is equivalent to a finite subsystem.
\end{theorem}
Note that the above result is only proved in \cite{Sela:hyp} for equations with no coefficients.  However, the more general result may be proved with exactly the same proof using the restricted Makanin--Razborov diagram rather than the (ordinary) Makanin--Razborov diagram.

Translating Theorem \ref{t:Noether} into the language of homomorphisms ensures that approximations will always exist.

\begin{theorem}\label{t:NoetherHom}
Suppose that $\Gamma$ is a torsion-free hyperbolic group.  For any sequence of epimorphisms of groups
\[
G_1\to G_2\to G_3\to\ldots
\]
the corresponding sequence of homomorphism varieties
\[
\Hom(G_1,\Gamma)\supset\Hom(G_2,\Gamma)\supset\Hom(G_3,\Gamma)\ldots
\]
eventually stabilizes.

In particular, for any finitely generated group $G$ there is a finitely presented group $\hat{G}$ along with an epimorphism $\eta \co \hat{G} \to G$ so that the injection $\eta^* \co \Hom(G,\Gamma) \to \Hom(\hat{G},\Gamma)$ is bijective.
\end{theorem}

If $G$ has a cyclic splitting, we can choose a $\Gamma$-approximation $\hat{G}$ that admits a splitting of the same form.   To do this, we give a recipe for building a $\Gamma$-approximation $\hat{G}$ for a group $G$, given a graph-of-groups decomposition $\Lambda$ of $G$ with cyclic edge groups, and $\Gamma$-approximations for the vertex groups of $\Lambda$.

\begin{lemma}\label{l:AFP approximation}
Suppose $G = A*_C B$ is a splitting of $G$ over a (possibly trivial) cyclic subgroup $C$.  Let $\eta_A:\hat{A}\to A$ and $\eta_B:\hat{B}\to B$ be $\Gamma$-approximations and let $\hat{C}_A\subseteq\hat{A}$ and $\hat{C}_B\subseteq\hat{B}$ be (possibly trivial) cyclic subgroups with the property that $\eta_A(\hat{C}_A)=\eta_B(\hat{C}_B)=C$.  Then
\[
\hat{G}=\hat{A}*_{\hat{C}_A=\hat{C}_B} \hat{B}
\]
is a $\Gamma$-approximation to $G$.
\end{lemma}
\begin{proof}
We can characterize $\Hom(G,\Gamma)$ as
\[
\Hom(G,\Gamma)=\{(\alpha,\beta)\in\Hom(A,\Gamma)\times\Hom(B,\Gamma)\mid \alpha|_{C}=\beta|_{C}\}.
\]
We can characterize $\Hom(\hat{G},\Gamma)$ similarly, and the result follows immediately.
\end{proof}

Likewise, we have a similar result for HNN-extensions.

\begin{lemma}\label{l:HNN approximation}
Suppose $G = A*_{C\sim C'}$ is a splitting of $G$, where the stable letter conjugates the infinite cyclic subgroups $C$ and $C'$.  Let $\eta_A:\hat{A}\to A$ be a $\Gamma$-approximation and let $\hat{C}\subseteq\hat{A}$ and $\hat{C}'\subseteq\hat{A}$ be cyclic subgroups with the property that $\eta_A(\hat{C})=C$ and $\eta_A(\hat{C}')=C'$.  Then
\[
\hat{G}=\hat{A}*_{\hat{C}\sim\hat{C}'}
\]
is a $\Gamma$-approximation to $G$.
\end{lemma}
\begin{proof}
We can characterize $\Hom(G,\Gamma)$ as
\[
\Hom(G,\Gamma)=\{(\alpha,\gamma)\in\Hom(A,\Gamma)\times\Gamma\mid \alpha|_{C}=\iota_\gamma\circ\alpha|_{C'}\}.
\]
where $\iota_\gamma\in\Aut(\Gamma)$ is conjugation by $\Gamma$.  Just as before, $\Hom(\hat{G},\Gamma)$ admits a similar characterization and the result follows.
\end{proof}

\section{If there are infinitely many} \label{s:Infinite}
\subsection{A characterization in terms of splittings of quotients}
In this section we discuss the cases when $\IHom$ and $\O$ are infinite.  The idea is to characterize the existence of infinitely many homomorphisms in terms of splittings of quotients of $H_\Sigma$.

\begin{proposition} \label{p:Inf2Limit}
If $\IHom$ is infinite then there is a strict restricted $\Gamma$-limit quotient $\eta \co H_\Sigma \to L$ so that infinitely many elements of $\IHom$ factor through $\eta$.  The quotient $L$ admits a free or essential cyclic splitting, and the coefficient subgroup is elliptic in this splitting.

If $\O$ is infinite then there is a strict $\Gamma$-limit quotient $\eta \co H_\Sigma \to L$ so that for infinitely many $o \in \O$, every representative of $o$ factors through $\eta$.  The quotient $L$ is either abelian or admits a free or essential cyclic splitting.
\end{proposition}
\begin{proof}
Consider the case with coefficients.  Assuming $\IHom$ is infinite, let $h_n \in \IHom$ be a sequence of pairwise distinct homomorphisms.  Passing to a subsequence, we may assume that $(h_n)$ is a stable sequence.  Let $L=H_\Sigma/\SK (h_n)$.  By \cite[Theorem 1.17]{Sela:hyp}, we may pass to a further subsequence and assume that every $h_n$ factors through the quotient map $\eta: H_\Sigma\to L$.  Since the $h_n$ are distinct, $L$ is a strict restricted $\Gamma$-limit group.  The quotient $L$ is non-abelian and so by Theorem \ref{t:CyclicSplit} admits a free or essential cyclic splitting in which the coefficient subgroup is elliptic.

In the coefficient-free case, the fact that  $\O$ is infinite implies that there exists a sequence of pairwise non-conjugate homomorphisms $h_n:H_\Sigma\to\Gamma$.  As in the previous paragraph, we can pass to a stable subsequence such that every $h_n$ factors through $\eta:H_\Sigma\to L=H_\Sigma/\SK (h_m)$.  Once again, Theorem \ref{t:CyclicSplit} ensures that if $L$ is non-abelian then $L$ splits as required.
\end{proof}

In fact, Proposition \ref{p:Inf2Limit} has a strong converse.  A simple version of this converse is encapsulated in the following two results.

\begin{proposition} \label{p:Limit2Inf-Coeff}
Suppose that $\Sigma$ and $\Lambda$ have coefficients.  Suppose $H_\Sigma$ admits a strict $\Gamma$-limit quotient $\eta : H_\Sigma \to L$ so that
\begin{enumerate}
\item some element of $\IHom$ factors through $\eta$; and
\item $L$ admits a free or essential cyclic splitting in which the coefficient subgroup is elliptic.
\end{enumerate}
Then $\IHom$ is infinite.
\end{proposition}

\begin{proposition} \label{p:Limit2Inf-NoCoeff}
Suppose that $\Sigma$ and $\Lambda$ are coefficient-free.  Suppose $H_\Sigma$ admits a strict $\Gamma$-limit  quotient $\eta : H_\Sigma \to L$ for which
\begin{enumerate}
\item a representative of some element of $\O$ factors through $\eta$; and
\item $L$ is abelian or admits a free or essential cyclic splitting.
\end{enumerate}
Then $\O$ is infinite.
\end{proposition}

We will prove stronger results than these, which are more easily applied in the algorithmic context (Propositions \ref{p:IHomisInfinite} and \ref{p:OisInfinite}).  However, for the moment it is useful to think about how one might go about proving these results.  The strategy is to pre-compose by Dehn twists, generalized Dehn twists or partial conjugations (see Definition \ref{d:Dehn twists etc}) to explicitly exhibit infinitely many elements of $\IHom$ or $\O$.  If we are considering only equations, this is straightforward.  In general, however, we need to ensure that the inequalities are preserved when twisting by non-inner automorphisms of $L$.

In order to do this, we use the following result.  It is well known in the case where $\Gamma$ is a free group, and is often called Baumslag's Lemma.  (Note that our proof can be substantially simplified in the free case.)

\begin{lemma}\label{l:Folklore lemma?}
Suppose $z, a_1,\ldots,a_n\in\Gamma$ and $[a_i,z]\neq 1$ for all $i$.  Then whenever $|m_i|$ are sufficiently large for $i< n$, the word
\[
w=a_1z^{m_1}\ldots a_nz^{m_n}
\]
does not commute with $z$.
\end{lemma}

The proof uses the $\delta$-hyperbolic geometry of the Cayley graph of $\Gamma$.  Let $d_\mathcal{A}$ be the word metric on the Cayley graph of $\Gamma$ with respect to $\mathcal{A}$ and suppose this metric is $\delta$-hyperbolic.  To prove Lemma \ref{l:Folklore lemma?} we need a version of the Margulis Lemma for hyperbolic groups.  The following is essentially contained in \cite{Paulin}.

\begin{lemma}\label{l:Margulis Lemma}
Let $a,b,c,d\in\Gamma$ and suppose that $c^a$ does not commute with $d^b$.  If $m$ and $n$ are integers, write $I_m=a[1,c^m]$ and $J_n=b[1,d^n]$ for choices of geodesic arcs in the Cayley graph of $\Gamma$.   Fix $\epsilon\geq 0$.  Then the diameter of
\[
I_m\cap N_\epsilon(J_n)
\]
(where $N_\epsilon(J_n)$ denotes the $\epsilon$-neighbourhood of $J_n$) is uniformly bounded, independently of $m$ and $n$.
\end{lemma}
\begin{proof}
Let $C_0$ be the number of elements in the ball of radius $8\delta+4\epsilon$.  Let $j\leq C_0$.

Let $I'\subseteq I_m$ be the arc of points at distance at least $|c|+2\delta$ from the endpoints of $I_m$.  Then for $x\in I'$, $d_{\mathcal{A}}((c^a)^i x,I_m)\leq 2\delta$.  Let $\pi$ be (a choice of) closest-point projection to $I_m$.  Then $\pi\circ (c^a)^i$ is at distance at most $2\delta$ from an orientation-preserving isometry $c':I'\to I_m$.

Similarly, let $J'\subseteq I_m$ be an arc of points at distance at least $|d^j|+2\delta+2\epsilon$ from the endpoints of $I_m$, and also within $\epsilon$ of $J_n$.  Then $\pi\circ (d^b)^j$ is at distance at most $2\delta+2\epsilon$ from an isometry $d'_j:J'\to I_m$.  Note that the translation length of $d'_j$ is bounded above by a linear function of $j$.

For a contradiction, suppose that the diameter of $I_m\cap N_\epsilon(J_n)$ is large enough that there exists $x\in I'\cap J'$ such that $c'(x)\in J'$ and $d'_j(x)\in I'$ for all $j\leq C_0$.  In particular, $d'_jc'(x)=c'd'_j(x)\in I_m$; we have that $d_{\mathcal{A}}((d^b)^jc^a(x),d'_jc'(x))\leq 4\delta+2\epsilon$ and, likewise, $d_{\mathcal{A}}(c^a(d^b)^j(x),c'd'_j(x))\leq 4\delta+2\epsilon$.  It follows that $d_{\mathcal{A}}((d^b)^{-j}(c^a)^{-1}(d^b)^jc^a(x),x)\leq 8\delta+4\epsilon$ and so, by the freeness of the action on the Cayley graph
\[
(d^b)^{-j}(c^a)^{-1}(d^b)^jc^a=(d^b)^{-j'}(c^a)^{-1}(d^b)^{j'}c^a
\]
for some distinct $j,j'\leq C_0$.  Using the fact that torsion-free hyperbolic groups are commutative transitive, it follows that $c^a$ commutes with $d^b$, as required.
\end{proof}

\begin{proof}[Proof of Lemma \ref{l:Folklore lemma?}]
We use the word-hyperbolic geometry of $\Gamma$.  Suppose the word metric $d_{\mathcal{A}}$ on $\Gamma$ corresponding to the generating set $\mathcal{A}$ is $\delta$-hyperbolic and denote the word length of an element $\gamma$ by $|\gamma|$.

The proof is by induction on $n$---note that the result is immediate for $n=1$.  If $w$ and $z$ commute then, by modifying $m_n$, we can assume that $w=1$.  Fix geodesic arcs $[1,z]$ and $[1,a_i]$ for each $i$. Write $w_i=a_1z^{m_1}\ldots z^{m_{i-1}}a_i$ so, by induction, we may assume that $w_i$ and $z$ do not commute for all $i$.  Let
\[
\alpha_i=w_{i-1}z^{m_{i-1}}[1,a_i]
\]
(where $w_0$ is taken to be $1$) and
\[
\zeta_i=w_i[1,z^{m_i}]
\]
so, if $w=1$ the concatenation $\alpha_1\cdot\zeta_1\cdot\alpha_2\cdots\alpha_n\cdot\zeta_n$ defines a $2n$-sided polygon in $\Gamma$.

Any such polygon is uniformly $2(n-1)\delta$-slim, so in particular $\zeta_1$ lies within the $2(n-1)\delta$-neighbourhood of the union of the $\alpha_i$ (over all $i$) and $\zeta_i$ (for $i>1$).  A segment of length at most $|a_i|+4(n-1)\delta$ can lie within the neighbourhood of $\alpha_i$, so if $A=\sum_i|a_i|$ then the lengths of segments of $\zeta_1$ covered by neighbourhoods of the other $\zeta_i$ sum to at least $|z^{m_1}|-A-4n(n-1)\delta$.  It follows easily from the fact that rectangles are $4\delta$-slim that, if two geodesic segments have endpoints that lie within $2(n-1)\delta$ of each other, then the entire segments lie within $2n\delta$ of each other.   Hence there exists a segment of $\zeta_1$ of length at least
\[
\frac{1}{n}\left(|z^{m_1}|-A-4n(n-1)\delta\right)
\]
that lies within a $2n\delta$-neighbourhood of some other $\zeta_i$.  Therefore by Lemma \ref{l:Margulis Lemma}, if $|m_1|$ is sufficiently large $a_1^{-1}w_i$ commutes with $z$, contradicting the inductive hypothesis.
\end{proof}

The following propositions are easy consequences.

\begin{proposition} \label{p:Dehn-nontriv}
Suppose that $a_1, \ldots , a_k, b_1, \ldots , b_k, c$ are elements of $\Gamma$ so that $a_1b_1 \cdots a_kb_k \ne 1$.
Then for all but finitely many $j$, the element \[  (c^{-j}a_1c^j)b_1 \cdots (c^{-j}a_kc^j)b_k     \] is non-trivial in $\Gamma$.
\end{proposition}

\begin{proposition} \label{p:Dehn-nontriv2}
Suppose that $a_1, \ldots , a_k, t,c$ are elements of $\Gamma$ so that $a_1t^{n_1} \cdots a_kt^{n_k} \neq 1$.  Then for all but finitely many $j$, the element \[  a_1(tc^j)^{n_1} \cdots a_k(tc^j)^{n_k}     \] is non-trivial in $\Gamma$.
\end{proposition}

\subsection{Approximate quotients}

In the algorithmic setting, the problem with applying Propositions \ref{p:Limit2Inf-Coeff} and \ref{p:Limit2Inf-NoCoeff} is that in general the limit quotient $L$ may not be finitely presented.  We deal with this by using $\Gamma$-approximations.

\begin{remark}\label{r:FP nearly-quotient}
Let $\pi:H_\Sigma\to L$ be any quotient of $H_\Sigma$ and let $\eta:\hat{L}\to L$ be a $\Gamma$-approximation to $L$.  By definition $\eta^*:\Hom(L,\Gamma)\to\Hom(\hat{L},\Gamma)$ is a bijection, and $\pi^*:\Hom(L,\Gamma)\to\Hom(H_\Sigma,\Gamma)$ is an injection.  So the composition
\[
\pi^*\circ(\eta^*)^{-1}:\Hom(\hat{L},\Gamma)\to\Hom(H_\Sigma,\Gamma)
\]
is an injection.

In the coefficient-free case, if $\hat{L}\cong\langle X \mid \mathcal{R}\rangle$ then the sentence
\[
\forall \underline{x}~\left(\mathcal{R}(\underline{x})\Longrightarrow\Sigma(\underline{x})\right)
\]
holds over $\Gamma$.  Likewise, in the case with coefficients, if $\hat{L}\cong\langle X,\mathcal{A}\mid \mathcal{R}\rangle$ then the sentence
\[
\forall \underline{x}~\left(\mathcal{R}(\underline{x},\underline{a})\Longrightarrow\Sigma(\underline{x},\underline{a})\right)
\]
holds in $\Gamma$.

Therefore, in order to show that $\IHom$ and $\O$ are infinite, we will show that there are infinitely many elements of $\mathrm{IHom}_\Gamma(\hat{L}, \Gamma ; \Lambda)$ and $\mathcal{O}(\hat{L}, \Lambda , \Gamma)$.
\end{remark}

\subsection{Dehn twists, generalized Dehn twists and partial conjugations}
\label{ss:Dehn,etc}

In this subsection we explain why the fact that the group $\hat{L}$ has a splitting is enough to deduce that $\IHomL$ and $\OL$ are infinite.

\begin{definition}
Suppose that $G$ is a finitely generated non-abelian group with a free or essential cyclic splitting $\Theta$.  A homomorphism $h : G \to \Gamma$ is called a {\em witness} with respect to $\Theta$ if
\begin{enumerate}
\item $h$ is injective on edge groups $C$;
\item $h(G)$ is non-abelian;
\item if $A$ is a non-abelian vertex group and the edge group $C$ is non-trivial then $h(C)$ is non-central in $h(A)$; and
\item if $\Theta$ is a free splitting then the image of each vertex group and each stable letter under $h$ is non-trivial.  Furthermore, if $A$ is a non-abelian
vertex group, $h(A)$ is non-abelian.
\end{enumerate}
\end{definition}

\begin{remark}\label{r:Limit groups have witnesses}
If $G$ is a non-abelian strict $\Gamma$-limit group then a witness always exists.
\end{remark}

A witness reflects the structure of the splitting $\Theta$ well enough that Dehn twists, generalized Dehn twists and partial conjugation give new (conjugacy classes of) homomorphisms to $\Gamma$.

\begin{definition}\label{d:Dehn twists etc}
There are sets of groups automorphisms associated to one-edge free or essential cyclic splittings.
\begin{enumerate}
\item If $G=A* B$ and $g\in A$ we can define the \emph{partial conjugation} $\alpha_g\in\Aut(H)$ such that $\alpha_g(a)=a$ for all $a\in A$ and $\alpha_g(b)=gbg^{-1}$ for all $b\in B$.
\item Suppose $G=A*_C B$ and $C$ is not central in $B$. If $C=\langle c\rangle$ then we can define the \emph{Dehn twist} $\delta_c\in\Aut(G)$ by $\delta_c(a)=a$ for $a\in A$ and $\delta_c(b)=c bc^{-1}$ for each $b\in B$.
\item Suppose $G=A*_C B$ and $B$ is abelian and torsion-free.  A \emph{generalized Dehn twist} is an automorphism $\delta'\in\Aut(G)$ that acts trivially on $A$ and that acts via a unimodular linear automorphism of $B$, fixing $C$.
\item Suppose $G=A*_C$ and let $t$ be the stable letter.  If $C=\langle c\rangle$ then we can define the \emph{Dehn twist} $\delta_c\in\Aut(G)$ by $\delta_c(a)=a$ for $a\in A$ and $\delta_c(t)=tc$.
\end{enumerate}
\end{definition}

\begin{proposition} \label{p:IHomisInfinite}
Suppose that $G$ is a group and that $\Theta$ is a non-trivial (one-edge) free or essential cyclic splitting of $G$.  Let $X\cup\mathcal{A}$ be a finite generating set for $G$, and assume that the coefficient subgroup $\langle\mathcal{A}\rangle$ is elliptic in $\Theta$.  Suppose furthermore that $\Lambda(\underline{x},\underline{a})\neq 1$ is a finite system of inequations, and that there is a homomorphism $h : G \to \Gamma$ so that
\begin{itemize}
\item[(i)] $h$ is a witness for $\Theta$; 
\item[(ii)] the restriction of $h$ to the coefficient subgroup $\langle\mathcal{A}\rangle$ surjects onto $\Gamma$; and
\item[(iii)] the inequations $\Lambda(h(\underline{x}),h(\underline{a}))\neq 1$ hold in $\Gamma$.
\end{itemize}
Then $\mathrm{IHom}_\Gamma(G,\Gamma;\Lambda)$ is infinite.
\end{proposition}
\begin{proof}
The idea is to twist $h$ by partial conjugations, Dehn twists or generalized Dehn twists to construct infinitely many elements of $\Hom(G,\Gamma)$ and then to use Propositions \ref{p:Dehn-nontriv} and \ref{p:Dehn-nontriv2} to conclude that all but finitely many of these elements satisfy the inequalities $\Lambda\neq 1$.  There are four cases to consider.

If $\Theta$ is a free splitting $G=A*B$ then we can assume that the coefficient subgroup is contained in $A$.  Note that there is an element $g\in A$ such that $h(g)$ does not centralize $h(B)$.  The sequence of homomorphisms $(h\circ\alpha_g^n)_{n\in\N}$ exhibits infinitely many elements of $\Hom_\Gamma(G,\Gamma)$ and, by Proposition \ref{p:Dehn-nontriv}, the inequalities $\Lambda(h\circ\alpha_g^n(\underline{x}),h(\underline{a}))\neq 1$ do not hold for at most finitely many $n$.  Therefore $\IHomG$ is infinite.

Suppose that $\Theta$ is an essential cyclic splitting of the form $G=A*_C B$ with $A$ and $B$ non-abelian and the coefficient subgroup is contained in $A$.  If $C=\langle c\rangle$ then $h(c)$ does not lie in the centre of $h(B)$ so, exactly as in the case of the free product, infinitely many of the homomorphisms $(h\circ\delta^n_c)_{n\in\N}$ are in $\IHomG$.

Suppose now that $\Theta$ is an essential cyclic splitting of the form $G=A*_C B$ with $A$ non-abelian and $B$ abelian, and the coefficient group is contained in $A$.  As the splitting is essential, $B$ is of rank at least 2.  Write $B = B_0 \oplus \langle t \rangle$ where $h(t) \ne 1$ and $C \subseteq B_0$.  Define the generalized Dehn twist
$\delta' \in \Aut(G)$ by $\delta'(a) = a$ for $a \in A$, $\delta'(b) = b$ for $b \in B_0$ and $\delta'(t) = tc$, where $c \in C \ssm \{ 1 \}$.  Note that since $h$ is a witness, $h \ne h \circ \delta'$.
The sequence of homomorphisms $(h\circ(\delta')^n)_{n\in\N}$ exhibits infinitely many elements of $\Hom_\Gamma(G,\Gamma)$ and, by Proposition \ref{p:Dehn-nontriv2}, the inequalities $\Lambda(h\circ(\delta')^n(\underline{x}),h(\underline{a}))\neq 1$ do not hold for at most finitely many $n$.  Therefore $\IHomG$ is infinite.

Suppose finally that $\Theta$ is a cyclic splitting of the form $G=A*_C$.  If $C=\langle c\rangle$ then the sequence of homomorphisms $(h\circ\delta_c^n)_{n\in\N}$ exhibits infinitely many elements of $\Hom_\Gamma(G,\Gamma)$.  As before, it follows from Proposition \ref{p:Dehn-nontriv2} that $\IHomG$ is infinite.
\end{proof}

We have a similar result for $\OG$.

\begin{proposition} \label{p:OisInfinite}
Suppose that $G$ is a non-abelian group and that $\Theta$ is a non-trivial (one-edge) free or essential cyclic splitting of $G$, and $\Theta$ has at least one non-abelian vertex group.  Let $X$ be a finite generating set for $G$.  Suppose furthermore that $\Lambda(\underline{x})\neq 1$ is a finite system of inequations, and that there is a homomorphism $h : G \to \Gamma$ so that
\begin{itemize}
\item[(i)] $h$ is a witness for $\Theta$; and
\item[(ii)] the inequations $\Lambda(h(\underline{x}))\neq 1$ hold in $\Gamma$.
\end{itemize}
Then $\mathcal O(G,\Lambda,\Gamma)$ is infinite.
\end{proposition}

The proof is identical to the proof of Proposition \ref{p:IHomisInfinite}.  After noting that the centralizer of any non-abelian subgroup of $\Gamma$
is trivial, it is straightforward to verify that the homomorphisms constructed are pairwise non-conjugate.

\begin{remark}\label{r:abelian case}
If $G$ is abelian and $\OG$ contains a non-trivial conjugacy class of homomorphisms then $\OG$ is infinite.  Similarly, if $G$ is a free product of
abelian groups and $\OG$ is non-empty then $\OG$ is infinite.
\end{remark}

\subsection{Checking for infinitely many solutions}

We have assembled all the tools we need to check that $\IHom$ and $\O$ are infinite.

\begin{definition}
A presentation for a group $G$ is said to \emph{exhibit a free splitting} $A*B$ if it is of the form
\[
\langle \mathcal{A},\mathcal{B}\mid\mathcal{R}_{\mathcal{A}},\mathcal{R}_{\mathcal{B}}\rangle
\]
where $\mathcal{R}_{\mathcal{A}}$ involves only generators in $\mathcal{A}$ and $\mathcal{R}_{\mathcal{B}}$ involves only generators in $\mathcal{B}$ (and $A\cong \langle \mathcal{A}\mid\mathcal{R}_{\mathcal{A}}\rangle$ and $B\cong \langle \mathcal{B}\mid\mathcal{R}_{\mathcal{B}}\rangle$).

Likewise, a presentation for a group $G$ is said to \emph{exhibit a cyclic splitting} $A*_CB$ if it is of the form
\[
\langle \mathcal{A},\mathcal{B},c\mid\mathcal{R}_{\mathcal{A}},\mathcal{R}_{\mathcal{B}}, c=w_{\mathcal{A}}, c=w_{\mathcal{B}}\rangle
\]
where $\mathcal{R}_{\mathcal{A}}$ involves only generators in $\mathcal{A}$, $\mathcal{R}_{\mathcal{B}}$ involves only generators in $\mathcal{B}$, $w_{\mathcal{A}}$ is a word in the generators $\mathcal{A}$ and $w_{\mathcal{B}}$ is a word in the generators $\mathcal{B}$.  (Here $A\cong \langle \mathcal{A}\mid\mathcal{R}_{\mathcal{A}}\rangle$, $B\cong \langle \mathcal{B}\mid\mathcal{R}_{\mathcal{B}}\rangle$ and $C=\langle c\rangle$.)

There is an obvious analogous definition of what it means for a presentation to exhibit an HNN-extension over a cyclic subgroup.
\end{definition}

In what follows, any vertex group is treated as non-abelian unless its relations include commutators of each pair of generators.

\begin{definition}
A presentation that exhibits a cyclic splitting is said to \emph{exhibit an essential cyclic splitting} if whenever a vertex group is visibly abelian it is not virtually cyclic.
\end{definition}

Note that it is easy to check algorithmically if a presentation exhibits a free or essential cyclic splitting, and to check whether a vertex group is visibly abelian.  Of course, one cannot in general decide if the edge group of a cyclic splitting is infinite, or indeed if a splitting is non-trivial. In our algorithms we will find witnesses that prove that our splittings are indeed as they appear.

\begin{theorem}\label{t:Checking if IHom is infinite}
Let $\Gamma$ be a torsion-free hyperbolic group and $\Sigma=1,\Lambda\neq 1$ a finite system of equations and inequations.  There exists a Turing machine that takes as input $\Sigma$ and $\Lambda$ and terminates if and only if the set of solutions $\IHom$ is infinite.
\end{theorem}
\begin{proof}
The Turing machine systematically enumerates group presentations of the form
\[
\hat{L}\cong\langle X,\mathcal{A}\mid\mathcal{R} \rangle
\]
and checks to see whether $\hat{L}$ satisfies the following conditions:
\begin{enumerate}
\item the presentation exhibits a free or essential cyclic splitting $\Theta$ of $\hat{L}$ in which the coefficient subgroup is elliptic;
\item  for each vertex group $V$ of $\Theta$ which is not visibly abelian there is a homomorphism $f \co \hat{L} \to \Gamma$ so that $f(V)$ is non-abelian;
\item there exists a witness $h$ for $\Theta$ that satisfies the inequations $\Lambda(h(\underline{x}),h(\underline{a}))\neq 1$---note that the existence of a witness can be checked by solving a finite system of equations and inequations over $\Gamma$; and
\item the sentence
\[
\forall \underline{x}~\left(\mathcal{R}(\underline{x},\underline{a})\Longrightarrow\Sigma(\underline{x},\underline{a})\right)
\]
holds in $\Gamma$.
\end{enumerate}
If such a group $\hat{L}$ is found then $\IHomL$ is infinite, by Proposition \ref{p:IHomisInfinite}.  But, since
\[
\curlyR(\underline{x},\underline{a})\Longrightarrow\Sigma(\underline{x},\underline{a})
\]
it follows that $\IHom$ is infinite by Remark \ref{r:FP nearly-quotient}.

It remains to show that such an $\hat{L}$ always exists if $\IHom$ is infinite.  By Proposition \ref{p:Inf2Limit}, $H_\Sigma$ has a limit group quotient $L$ with a free or essential cyclic splitting in which the coefficient group is elliptic.  By Lemmas  \ref{l:AFP approximation} and \ref{l:HNN approximation}, there exists a $\Gamma$-approximation
\[
\hat{\eta}:\hat{L}\to L
\]
with a splitting of the same type, in which the coefficient subgroup is elliptic.  Since $L$ is a limit group there exists a witness $h:L\to\Gamma$, and the homomorphism $h\circ\hat{\eta}$ is easily seen to be a witness for the splitting of $\hat{L}$. By Remark \ref{r:FP nearly-quotient}, if $\hat{L}$ has presentation $\langle X,\mathcal{A}\mid\mathcal{R} \rangle$ then the identity
\[
\curlyR(\underline{x},\underline{a})\Longrightarrow\Sigma(\underline{x},\underline{a})
\]
holds in $\Gamma$.  Therefore, there does exist such a group $\hat{L}$ and so a systematic search will eventually find it and the proof that it has the required properties.
\end{proof}

Likewise, we have the following theorem.

\begin{theorem}\label{t:Checking if O is infinite}
Let $\Gamma$ be a torsion-free hyperbolic group and $\Sigma=1,\Lambda\neq 1$ a finite, coefficient-free system of equations and inequations.  There exists a Turing machine that takes as input $\Sigma$ and $\Lambda$ and terminates if and only if the set of conjugacy classes of solutions $\O$ is infinite.
\end{theorem}
\begin{proof}
The proof is very similar to the proof of Theorem \ref{t:Checking if IHom is infinite}.  The Turing machine systematically enumerates finite presentations of the form
\[
\hat{L}\cong\langle X\mid\mathcal{R} \rangle.
\]
For each presentation, the Turing machine checks to see whether $\hat{L}$ satisfies the following conditions.
\begin{enumerate}
\item The presentation should exhibit a free or essential cyclic splitting $\Theta$ of $\hat{L}$, or exhibit $\hat{L}$ as an abelian group.
\begin{enumerate}
\item In the split case,  vertex groups which are not visibly abelian should be seen to be non-abelian via a homomorphism to $\Gamma$.
Furthermore,
 there should exist a witness $h$ for $\Theta$ that satisfies the inequations $\Lambda(h(\underline{x}))\neq 1$---again the existence of a witness can be checked by solving a finite system of equations and inequations over $\Gamma$.
\item In the abelian case, there should exist a non-trivial element of $\OL$.
\end{enumerate}
\item Finally, the sentence
\[
\forall \underline{x}~\left(\mathcal{R}(\underline{x})\Longrightarrow\Sigma(\underline{x})\right)
\]
must hold in $\Gamma$.
\end{enumerate}
If such a group $\hat{L}$ is found then $\O$ is infinite, by Proposition \ref{p:OisInfinite} or by Remark \ref{r:abelian case}.  But, since
\[
\curlyR(\underline{x})\Longrightarrow\Sigma(\underline{x})
\]
it follows that $\O$ is infinite.  That there exists such a limit group quotient if $\O$ is infinite follows exactly as in the proof of Theorem \ref{t:Checking if IHom is infinite}.
\end{proof}

\section{If there are finitely many} \label{s:Finite}

In this section we describe the other halves of the algorithms in the statements of Theorem \ref{t:FinSol} and Theorem \ref{t:FinConjSol}. We start with the case when $\IHom$ is finite, as this is much easier.

\subsection{Finitely many solutions}

\begin{theorem}\label{t:Finitely many solutions}
Let $\Gamma$ be a torsion-free hyperbolic group. There is an algorithm which takes as input a finite system of equations and inequations (with coefficients) over $\Gamma$ and terminates if and only if the system has only finitely many solutions.

If it terminates, it provides a list of the solutions.
\end{theorem}
\begin{proof}
Let $\Sigma = 1$, $\Lambda \ne 1$ be the finite system.

According to Sela \cite[Theorem 7.12]{Sela:hyp} (or alternatively \cite[Theorem 0.1]{Dah_eq}), we can determine if there is a solution or not.  If there is not, we terminate with the empty list.

If there is a solution, we search for it systematically.  Since we know that there is a solution, we will eventually find it.  The fact that there is a solution which does not appear on a given finite list is easily encoded into a finite system of equations and inequations, and so we can determine if there is a new solution.

Continuing in this manner, we will clearly discover if there are only finitely many solutions to $\Sigma =1, \Lambda \ne 1$, and in case we discover this we will find a list of all of the solutions.
\end{proof}

Combining Theorem \ref{t:Finitely many conjugacy classes of solutions} with Theorem \ref{t:Checking if IHom is infinite} completes the proof of Theorem \ref{t:FinSol}.

\subsection{Finitely many conjugacy classes of solutions}

We now turn to the question of recognizing when $\O$ is finite.  This is substantially more difficult than recognizing when $\IHom$
is finite.  In spirit, we want to find solutions that are shortest in their conjugacy classes.  This is in general too much to ask, but the results of \cite{DG} provide us with a slightly weakened notion that is an adequate substitute.  We refer the reader to \cite[Section 4]{DG} for notation and terminology.  Note that we are not concerned with primary peripheral structures in this paper, and we can take them to be empty.

\begin{definition}
Let $H$ be a group with fixed generating set $X$.  A homomorphism $\psi:H\to\Gamma$ is \emph{compatible} if it is injective on the ball of radius $8$ in $H$.
\end{definition}

Compatible homomorphisms are problematic for two reasons.  First, a conjugacy class $[\phi]\in\O$ may not have a compatible representative.  Secondly, we may not be able to calculate the ball of radius $8$, as we have no solution to the word problem in $H_\Sigma$.  We will circumvent these by supplementing $\Sigma$ and $\Lambda$ in order to determine the ball of radius $8$.

\begin{definition}
Let $k$ be a natural number, and let $\mathcal{B}_{F(X)}(k)$ be the ball of radius $k$ in the free group $F(X)$.  A system of equations and inequations $\Sigma=1$, $\Lambda\neq 1$ \emph{forces the ball of radius $k$} if there is a subset $S\subseteq \mathcal{B}_{F(X)}(k)$ such that $S\subseteq \Sigma$ and $\mathcal{B}_{F(X)}(k)\ssm S\subseteq\Lambda$.
\end{definition}

There are finitely many subsets $S$ of the ball of radius $8$ in $F(X)$.  For such an $S$, let $\Sigma_S=\Sigma\cup S$ and let $\Lambda_S=\Lambda\cup(\mathcal{B}_{F(X)}(8)\ssm S)$.   For any solution $\phi$ to $\Sigma=1,\Lambda\neq 1$ there is a unique $S\subseteq\mathcal{B}_{F(X)}(8)$ such that $\phi$ satisfies $\Sigma_S=1$ and $\Lambda_S\neq 1$. Therefore we have the following remark.

\begin{remark}\label{r:Compatible decomposition}
For $S\subseteq\mathcal{B}_{F(X)}(8)$, let $\eta_S:H_\Sigma\to H_{\Sigma_S}$ be the obvious epimorphism.  The system $\Sigma_S=1,\Lambda_S\neq 1$ forces the ball of radius $8$.  Furthermore,
\[
\O=\coprod_{S\subseteq\mathcal{B}_{F(X)}(8)} \eta_S^*~\mathcal{O}(H_{\Sigma_S},\Gamma;\Lambda_S)
\]
Moreover, any representative $\phi$ of a conjugacy class $[\phi]\in\mathcal{O}(H_{\Sigma_S},\Gamma;\Lambda_S)$ is a compatible homomorphism $H_{\Sigma_S}\to\Gamma$.
\end{remark}

Note that $\O$ is finite if and only if $\mathcal{O}(H_{\Sigma_S},\Gamma;\Lambda_S)$ is finite for each $S\subseteq\mathcal{B}_{F(X)}(8)$.

The property $\Omega$ defined in \cite[Remark 4.8]{DG} provides an adequate means of recognizing short homomorphisms in conjugacy classes.  The property $\Omega$ applies to {\em acceptable lifts} of a 
homomorphism.  See \cite{Dah_eq} for the definition of an acceptable lift and also \cite{DG} for more details.  

\begin{definition}
Fix a pair of elements $a,b\in\Gamma$ that do not commute.  As $\Gamma$ is non-elementary we can always find such a pair.  A homomorphism $\psi:H_\Sigma\to\Gamma$ is \emph{fairly short} if it is compatible and has an acceptable lift $\tilde{\psi}: \mathcal{B}_{H_\Sigma}(2)\to F$ that satisfies $\Omega$.  The homomorphism $\psi$ is \emph{very short} if it is compatible and every acceptable lift $\tilde{\psi}: \mathcal{B}_{H_\Sigma}(2)\to F$ satisfies $\Omega$.
\end{definition}

\begin{remark}
Here are some important observations.
\begin{enumerate}
\item Note that $\Omega$, as defined in \cite[Remark 4.9]{DG}, depends on the choice of $a$ and $b$.
\item The condition $\Omega$ is a boolean combination of \emph{normalized rational constraints} on the lift $\tilde{\psi}$.
\end{enumerate}
\end{remark}

By \cite[Lemma 4.7 and Proposition 4.8]{DG}, we obtain the following.

\begin{lemma}\label{l:Finitely many fairly short homs}
Any conjugacy class of compatible homomorphisms $H_\Sigma\to\Gamma$ contains at least one very short homomorphism, and at most finitely many fairly short homomorphisms.
\end{lemma}

Fairly short and very short homomorphism are useful because, as they are characterized by a boolean combination of normalized rational constraints on acceptable lifts, we can follow \cite{Dah_eq} in exploiting the work of \cite{DGH} to identify all the very short solutions to a system of equations and inequations.  In this terminology, we can deduce the following from \cite[Proposition 1.5 ]{Dah_eq}.

\begin{proposition}\label{p:Finding short solutions}
There exists an algorithm that, given a finite presentation for a torsion-free hyperbolic group $\Gamma$ and a system of equations $\Sigma =1$ and inequations $\Lambda\neq 1$ with coefficients in $\Gamma$ always terminates.  If it answers ``yes" then there is a fairly short solution to $\Sigma=1, \Lambda \ne1$.  If it answers ``no" then there is no
very short solution.
\end{proposition}

This is what we need in order to test whether there are only finitely many conjugacy classes of solutions.

\begin{theorem}\label{t:Finitely many conjugacy classes of solutions}
Let $\Gamma$ be a torsion-free hyperbolic group.  There is an algorithm that takes as input a finite system of equations $\Sigma=1$ and inequations $\Lambda\neq 1$ (without coefficients) and terminates if and only if there are only finitely many conjugacy classes of solutions to the system.  In case the algorithm terminates, it provides a list consisting of exactly one representative of each conjugacy class of solutions.
\end{theorem}
\begin{proof}
In fact, we will describe an algorithm that provides a list containing at least one, and at most finitely many, representatives of each conjugacy class of solutions to $\Sigma=1,\Lambda\neq 1$.  It is easy to pare down such a list so that it contains precisely one representative of each conjugacy class of solutions, using a solution to the simultaneous conjugacy problem in $\Gamma$.

For each subset $S\subseteq\mathcal{B}_{F(X)}(8)$ in turn, consider the system $\Sigma_S=1,\Lambda_S\neq 1$.  By Remark \ref{r:Compatible decomposition}, $\O$ is infinite if and only if each of the corresponding $\mathcal{O}(H_{\Sigma_S},\Gamma;\Lambda_S)$ are infinite.  We may therefore assume that $\Sigma=1,\Lambda\neq 1$ forces the ball of radius 8.  In particular, all solutions are compatible.

Let us now describe the procedure for producing the list.  Let $l$ be our list of solutions, which initially is empty.  Add to $\Sigma=1,\Lambda\neq 1$ the stipulation that no solution should be equal to a solution on the list $l$, creating the new system of equations and inequations $\Sigma=1,\Lambda'\neq 1$.  Since we are only adding inequations, we have not changed the notions of very short and fairly short.  Now apply the algorithm from Proposition \ref{p:Finding short solutions} to the system $\Sigma=1,\Lambda'\neq 1$.  
If the algorithm from Proposition \ref{p:Finding short solutions} says ``no" then there is no very short solution to $\Sigma=1$ and $\Lambda'\ne 1$, and the procedure
terminates with the list $l$.    If the algorithm says ``yes" then there is a fairly short solution.  A naive search will eventually find such a solution.  It cannot be on the list $l$, so add it to $l$ and repeat the process.

By Lemma \ref{l:Finitely many fairly short homs}, every conjugacy class of solutions has a very short representative.  It follows that, if $\O$ is infinite then the algorithm will not terminate.  On the other hand Lemma \ref{l:Finitely many fairly short homs} states that every conjugacy class contains only finitely many fairly short representatives, so if $\O$ is finite then the algorithm will eventually terminate and output the list $l$.  
\end{proof}

Combining Theorem \ref{t:Finitely many conjugacy classes of solutions} with Theorem \ref{t:Checking if O is infinite} completes the proof of Theorem \ref{t:FinConjSol}.

\section{Immutable subgroups} \label{s:Immutable}

\subsection{Enumerating immutable subgroups}

\begin{definition}
Let $G$ be a group.  A finitely generated subgroup $H$ of $G$ is called {\em immutable} if there are finitely many injections $\phi_1, \ldots \phi_N : H \to G$ so that any injection $\phi : H \to G$ is conjugate to one of the $\phi_i$.
\end{definition}

We will see that all $\Gamma$-limit groups are either strict or immutable subgroups of $\Gamma$.  Therefore, in order to understand $\Gamma$-limit groups, or even just subgroups of $\Gamma$, we need to study immutable subgroups.  This is one of the big differences between $\Gamma$-limit groups where $\Gamma$ is a torsion-free hyperbolic group and ordinary limit groups defined over a free group.  Clearly, there are no non-trivial immutable subgroups of a free group.

\begin{lemma} \label{l:Splits}
A finitely generated subgroup of $\Gamma$ is immutable if and only if it does not admit a non-trivial free splitting or an essential splitting over $\Z$.
\end{lemma}
\begin{proof}
Let $H$ be a finitely generated subgroup of $\Gamma$.  If $H \cong \Z$ then $H$ is not immutable, since $\Z$ splits over $1$.

If $H$ admits a non-trivial free or essential cyclic splitting then, taking the inclusion $h:H\hookrightarrow \Gamma$ as a witness, it follows immediately by Proposition \ref{p:OisInfinite} that there are infinitely many conjugacy classes of homomorphisms $H\to\Gamma$.  As each is obtained by pre-composing $h$ with an automorphism they are all injective, and so $H$ is not immutable.

Suppose, then, that $H$ admits infinitely many conjugacy classes of embeddings into $\Gamma$.  We can therefore choose a sequence of pairwise non-conjugate embeddings $\{\phi_i:H\to \Gamma\}$.  After passing to a subsequence, the sequence $\phi_i$ can be taken to be stable.  As $H=H/\SK(\phi_i)$, this realizes $H$ as a \emph{strict} $\Gamma$-limit group.  It now follows from Theorem \ref{t:CyclicSplit} that $H$ splits essentially over a cyclic subgroup.
\end{proof}

\begin{remark}\label{r:Immutable vs Strict}
The proof of Lemma \ref{l:Splits} shows that any $\Gamma$-limit group is \emph{either} strict \emph{or} an immutable subgroup of $\Gamma$, and the two cases are mutually exclusive.
\end{remark}

\begin{example}\label{e:Hyperbolic manifolds are immutable}
Let $\Gamma$ be the fundamental group of a closed hyperbolic $n$-manifold where $n \ge 3$.  Then $\Gamma$ admits no non-trivial abelian splitting (see, for example, \cite[Theorem 1.6(i)]{Bel}) and so, by Lemma \ref{l:Splits}, is an immutable subgroup of itself.
\end{example}

In \cite{OW}, a version of the Rips construction is given in which the kernel has Kazhdan's property $(T)$, and in particular is immutable.  This provides many further examples of finitely generated (but not finitely presentable) immutable subgroups.

The following lemma provides a slightly weaker characterization of immutable subgroups.  The advantage is that it can be recognized algorithmically.

\begin{lemma} \label{l:FinCondition}
A finitely generated subgroup $H$ of $\Gamma$ is immutable if and only if there are finitely many homomorphisms $\rho_1, \ldots , \rho_K : H \to \Gamma$ and an integer $D \ge 1$, so that for any homomorphism $\rho : H \to \Gamma$, either
\begin{enumerate}
\item $\rho$ is conjugate to one of the $\rho_i$; or
\item $\rho$ is not injective on the ball of radius $D$ about $1$ in $H$.
\end{enumerate}
\end{lemma}
\begin{proof}
If $H$ is not immutable then it admits infinitely many conjugacy classes of embeddings, and so clearly does not satisfy the dichotomy in the statement of the lemma.

On the other hand, suppose that $H$ does not satisfy the statement of the lemma.  Then there is a sequence $\{ \phi_i : H \to \Gamma \}$ of non-conjugate homomorphisms so that $\phi_i$ is injective on the ball of radius $i$ about $1$ in $H$.  We may now pass to a stable subsequence with trivial stable kernel and deduce that $H$ is a strict $\Gamma$-limit group and therefore not immutable.
\end{proof}

Note that we do not insist that the homomorphisms $\rho_i$ are injective.  However, since $H$ is assumed to be a subgroup of $\Gamma$, necessarily at least one of the $\rho_i$ is an embedding.

\begin{example}\label{e:Finitely many conjugacy classes of solutions}
Fix a generating set $X$ for $H$ an immutable subgroup of $\Gamma$ (e.g. Example \ref{e:Hyperbolic manifolds are immutable}).  The property that a map from $X$ into $\Gamma$ extends to a homomorphism $h : H \to \Gamma$ that is injective on the ball of radius $D$ in the Cayley graph of $H$ with respect to $X$ is easily encoded as a system of (possibly infinitely many) equations and finitely many inequations over $\Gamma$.   By Lemma \ref{l:FinCondition}, there are only finitely many conjugacy classes of solutions.   By Theorem \ref{t:Noether}, we have an example of a finite system of equations and inequations with only finitely many conjugacy classes of solutions in $\Gamma$.

Note that if $H$ happens to be finitely presentable (as in Example \ref{e:Hyperbolic manifolds are immutable}), it is not necessary to appeal to Theorem \ref{t:Noether}.
\end{example}

One of the problems with finitely generated subgroups of torsion-free hyperbolic groups is that they need not be finitely presentable.  There are also torsion-free hyperbolic groups for which it is undecidable which finite subsets generate finitely presentable subgroups (see \cite[Theorem 1]{Bridson-Wise}).  We do not know the answer to the following question.

\begin{question}
Let $\Gamma$ be a torsion-free hyperbolic group.  Is there an algorithm that takes as input finite subsets $\mathcal A$ of $\Gamma$ so that $\langle \mathcal A \rangle$ is finitely presented and finds a finite presentation for $\langle \mathcal A \rangle$?
\end{question}

We now have the tools to prove that we can recognize algorithmically when a finite set generates an immutable subgroup (although not when it does not!).  Theorem \ref{t:ListImm} is an immediate corollary.

\begin{proposition}\label{p:Checking immutability}
There exists a Turing machine that takes as input a finite subset $\mathcal A\subseteq\Gamma$ and terminates if and only if the subgroup $H=\langle \mathcal A\rangle$ is immutable.
\end{proposition}
\begin{proof}
Our approach is to construct a sequence of epimorphisms $\hat{H}_D\to H$ (for any $D\in\N$) and to test the criterion of Lemma \ref{l:FinCondition} on each $\hat{H}_D$.

Let $\hat{H}_D$ be the group generated by $S$ with relations equal to all the loops in the ball of radius $D$ in $H$.  So there is an epimorphism $\eta_D:\hat{H}_D\to H$ that is an isometry between the balls of radius $D$.  Homomorphisms $\hat{H}_D\to\Gamma$ that are injective on the ball of radius $D$ in $\hat{H}_D$ can be characterized precisely by a system of equations and inequations over $\Gamma$---the equations are the relations of $\hat{H}_D$ and the inequations separate the points in the ball of radius $D$.  (Although we do not have a solution to the word problem in $\hat{H}_D$, we can still construct its ball of radius $D$ since it is, by definition, isometric to the ball of radius $D$ in $H$.)  Using the algorithm of Theorem \ref{t:FinConjSol} we can now determine whether this system of equations and inequations has finitely many conjugacy classes of solutions.  If it does then it follows immediately that $H$ is immutable.

Checking this for each $D$ in turn, we have a process that, if it terminates, confirms that $H$ is immutable.  It remains to show that the process terminates whenever $H$ is immutable.  By Theorem \ref{t:NoetherHom}, $\Hom(\hat{H}_D,\Gamma)=\Hom(H,\Gamma)$ for all sufficiently large $D$---that is, in our terminology, $\hat{H}_D$ is a $\Gamma$-approximation for $H$.  For large enough $D$ we will also have that $D$ satisfies the conclusions of Lemma \ref{l:FinCondition}.  For such a $D$, there are only finitely many conjugacy classes of homomorphisms $\hat{H}_D\to\Gamma$ that are injective on the ball of radius $D$.  Therefore, if $H$ is immutable then the process will terminate.
\end{proof}

\subsection{Questions}

We have proved that the set of immutable subgroups of $\Gamma$ is recursively enumerable.

\begin{question}
Is the set of immutable subgroups of $\Gamma$ recursive?
\end{question}

Here a class of groups is {\em recursive} if there is an algorithm that enumerates the class of groups so that the list contains a single representative of each isomorphism class.   This is essentially the same as asking whether the isomorphism problem is solvable for immutable subgroups. One would expect the answer to be negative, as the isomorphism problem for finitely generated subgroups of hyperbolic groups is unsolvable.

Applying the construction of \cite{OW} to a group with unsolvable word problem, we obtain the following.

\begin{proposition}\label{p:Immutable Membership}
There are torsion-free hyperbolic groups $\Gamma$ for which the membership problem for finitely generated immutable subgroups is unsolvable.
\end{proposition}

However, no immutable subgroup constructed in this way can be finitely presentable.

\begin{question}\label{q:Membership problem}
Is the membership problem for finitely presentable immutable subgroups solvable?
\end{question}

Question \ref{q:Membership problem} appears to be related to the isomorphism problem for finitely presented subgroups of torsion-free hyperbolic groups.

\begin{question}\label{q:Recognizing immutability}
Let $\Gamma$ be a (torsion-free) hyperbolic group.  Is there an algorithm that takes as input a finite subset of $\Gamma$ and decides whether or
not the subgroup generated is immutable?
\end{question}

Decision problems about finitely generated subgroups of hyperbolic groups are rarely solvable (see, for instance, \cite{Rips}, \cite{Bridson-Wise}).  Nevertheless, the results of this paper reduce Question \ref{q:Recognizing immutability} to determining if a subgroup of $\Gamma$ splits freely or essentially over $\Z$.

\bibliography{temp}

\begin{thebibliography}{99}
\bibitem{Bel} I. Belegradek, Aspherical manifolds with relatively
hyperbolic fundamental groups, \textit{IJAC} {\bf 18} (2008), 97--110.

\bibitem{Bestvina} M. Bestvina, Degenerations of the hyperbolic space, \textit{Duke Math. J.} {\bf 56} (1988), 143--161.

\bibitem{Bestvinahandbook} M. Bestvina, $\R$-trees in topology, geometry, and group theory, \textit{Handbook of geometric topology}, 55--91, North
Holland, 2002.

\bibitem{BFSela} M. Bestvina and M. Feighn, Notes on Sela's work: Limit groups
and Makanin-Razborov diagrams, preprint.

\bibitem{Bridson-Haefliger}
Martin~R. Bridson and Andr\'{e} Haefliger.
\newblock {\em Metric spaces of non-positive curvature}, volume 319 of {\em
  Grundlehren der mathematischen {W}issenschaften}.
\newblock Springer, 1999.

\bibitem{BS} M.R. Bridson and G.A. Swarup, On Hausdorff-Gromov convergence and a theorem of Paulin, \textit{Enseign. Math.} {\bf 40} (1994), 267--289.


\bibitem{Bridson-Wise} M. Bridson and D. Wise, Malnormality is undecidable in
hyperbolic groups, \textit{Israel J. Math.} {\bf 124} (2001), 313--316.

\bibitem{Dah_eq} F. Dahmani, Existential questions in (relatively) hyperbolic groups, \textit{Israel J. Math.}, to appear.

\bibitem{DG} F. Dahmani and D. Groves, The Isomorphism Problem for toral
relatively hyperbolic groups, \textit{Publ. Math. IHES} {\bf 107} (2008), 211--290..

\bibitem{DGH} V. Diekert, C. Guti\'errez, C. Hagenah,  The existential theory of equations with rational constraints in free groups is PSPACE-complete.  STACS 2001 (Dresden),  170--182, Lecture Notes in Comput. Sci., 2010, Springer, 2001.

\bibitem{MR-RH} D. Groves, Limit groups for relatively hyperbolic groups, II:
Makanin-Razborov diagrams, \textit{Geom. and Topol.} {\bf 9} (2005), 2319--2358.

\bibitem{GW1} D. Groves and H. Wilton, Enumerating limit groups, \textit{Groups, Geom. and Dynamics} {\bf 3} (2009), 389--399.

\bibitem{GW3} D. Groves and H. Wilton, Limit groups over hyperbolic groups, in preparation.

\bibitem{Guirardel} Vincent Guirardel.
\newblock Actions of finitely generated groups on $\mathbb{R}$-trees.
\newblock {\textit{Ann. Inst. Fourier}}, {\bf 58} (2008), 159--211.

\bibitem{KM} O. Kharlampovich and A. Miasnikov, Irreducible affine varieties
over a free group II, \textit{Journal of Algebra} {\bf 200} (1998), 517--570.

\bibitem{KM:Final} O. Kharlampovich and A. Miasnikov, Elementary theory of
nonabelian free groups, \textit{J. Algebra} {\bf 302} (2006), 451--552.

\bibitem{Makanin} G.S. Makanin, Equations in a free group. (Russian) \textit{Izv.
Akad. Nauk SSSR Ser. Nat.} {\bf 46} (1982), 1199--1273; translation in
\textit{Math. USSR-Izv.} {\bf 21} (1983), 546--582.

\bibitem{Makanin2} G.S. Makanin, Decidability of the universal and positive theories
of a free group. (Russian) \textit{Izv. Akad. Nauk SSSR Ser. Mat.} {\bf 48} (1984), 735--749; translation in \textit{Math. USSR-Izv.} {\bf 25} (1985), 75--88.

\bibitem{OW} Y. Ollivier and D. Wise, Kazhdan groups with infinite outer automorphism group, \textit{Trans. AMS} {\bf 360} (2007), 1959--1976.

\bibitem{Paulin} Fr{\'e}d{\'e}ric Paulin.
\newblock {\em Outer automorphisms of hyperbolic groups and small actions on
  $\mathbb R$-trees}, volume~19 of {\em Math. Sci. Res. Inst. Publ.}, pages
  331--343.
\newblock Springer, New York, 1991.

\bibitem{Razborov} A. Razborov, On systems of equations in a free group,
\textit{Math. USSR Izvestiya} {\bf 25} (1985), 115--162.

\bibitem{Razborov:thesis} A. Razborov, On systems of equations in a free group,
Ph.D. thesis, \textit{Steklov Math. Institute} (1987).

\bibitem{Rips} E. Rips, Subgroups of small cancellation groups, \textit{Bull LMS} {\bf 14} (1982), 45--47.

\bibitem{RS} E. Rips and Z. Sela, Canonical representatives and equations in hyperbolic groups.  \textit{Invent. Math.}  {\bf 120}  (1995),  no. 3, 489--512.

\bibitem{Sela:Inventiones} Z.~Sela.
\newblock Acylindrical accessibility for groups.
\newblock \textit{Invent. Math.}, {\bf 129} (1997) 527--565.

\bibitem{Sela:IHES} Z. Sela, Diophantine geometry over groups, I: Makanin-Razborov
diagrams, \textit{Publ. Math. IHES} {\bf 93} (2001), 31--105.

\bibitem{Sela:last} Z. Sela, Diophantine geometry over groups. VI. The elementary
theory of a free group, \textit{GAFA} {\bf 16} (2006), 707--730.

\bibitem{Sela:hyp} Z. Sela, Diophantine geometry over groups, VII: The elementary
theory of hyperbolic groups, \textit{Proc. LMS} {\bf 99} (2009), 217--273..






\end{thebibliography}

\end{document}